\documentclass[a4paper, 12pt]{article}
 \usepackage[utf8]{inputenc} 
\usepackage{makeidx}
\usepackage{enumerate, theorem}
\usepackage{amsmath, amsfonts, amssymb}
\usepackage[height=22.5cm, width=15cm]{geometry}
\usepackage[all]{xy}
\usepackage{mathrsfs, yfonts}
\usepackage{graphicx}
\DeclareMathOperator{\C}{\mathbb{C}}

\newcommand{\parag}[1]{\paragraph{\sc{#1.}} }

\newtheorem{thm}{Theorem}[subsection]
\newtheorem{defn}[thm]{Definition}

\newtheorem{prop}[thm]{Proposition}

\begin{document}

\title{Some comments on my study of period-integrals}

\date{19 Mars 2024}

\author{Daniel Barlet\footnote{Barlet Daniel, Institut Elie Cartan UMR 7502  \newline
Universit\'e de Lorraine, CNRS,  Institut Universitaire de France, \newline
BP 239 - F - 54506 Vandoeuvre-l\`es-Nancy Cedex.France. \newline
e-mail : daniel.barlet@univ-lorraine.fr}.}

\maketitle

\parag{Abstract} This text is a presentation (without proofs)  of some of my  recent results  on the singular terms of asymptotic expansions of period-integrals using
(a,b)-modules. I try to explain why this simple algebraic structure is interesting and useful.

\parag{AMS classification}  32 S 25; 32 S 40 ; 34 E 05

\parag{Key words} Asymptotic Expansions ; Period-integral;  Convergent (a,b)-Module; Hermitian Period;  Geometric (convergent) (a,b)-Modules; (convergent) Fresco; Bernstein Polynomial ; Higher Order  Bernstein Polynomials. \\

\section{Introduction}

The aim of this text is devoted to give some comments on my work since more than thirty years on the  local study of the singularity of a holomorphic function in a complex manifold.\\
As my point of view seems not very popular and looks at least strange and often few interesting for many specialists of this subject, I shall try to explain why it is pertinent and 
interesting\footnote{But the notion of "interesting result" is not so easy to define in an objective way. Only long range vision may give a serious answer to such a question.}. \\
In particular I shall try to explain how the use of the structure of (a,b)-module may be seen as a tentative of having a "local Hodge Theory" even if, at first glance, such a notion looks
meaning less. Note that I do not pretend that this point of view is an alternative to the use of desingularization and usual (global) Hodge Theory (or Mixte Hodge modules \cite{[S.90]}), but I think
that it may be a complementary approach and can give interesting and more precise  informations when we are interested in a specific period integral.\\

\section{The context}

Consider a holomorphic (non constant) germ $f :  (\C^{n+1}, 0) \to (\C, 0)$ with an isolated singularity for the eigenvalue $r :=exp(-2i\pi\alpha)$ of the monodromy (we assume
 that $\alpha$ is in $ ]0, 1]\cap \mathbb{Q}$). Then, for a given homology class $\gamma \in H_n(F_0, \C)_r$, where $F_0$ is the Milnor fiber of $f$ at $0$ and where $H_n(F_0, \C)_r$ is the
 generalized eigenspace for the eigenvalue $r$ of the monodromy of $f$ acting on $H_n(F_0, \C)$, we defined the period-integral
 \begin{equation}
 \Phi_{\omega, \gamma}(s) := \int_{\gamma(s)} \omega/df 
 \end{equation}
 where $\omega$ is any germ in $\Omega_0^{n+1}$, where $s$ is in $H$ the universal cover of a small punctured disc $D^*$ with center $0$ in $\C$ and where $(\gamma(s))_{s \in H}$ is the horizontal family
  of  compact $n$-cycles in the fibers of $f$ taking the value $\gamma$ at the base point $s_0$ in $H$.\\
  Such multivalued holomorphic function admits an  asymptotic expansion at the origin (convergent in sectors) for $\alpha, \gamma, \omega$ given,  which looks like
  \begin{equation}
  \Phi_{\omega, \gamma}(s) \simeq \sum_{j \in [0, n], m \in \mathbb{N}} \varphi_{m, j}(\gamma) s^{\alpha+ m-1}(Log \, s)^j
  \end{equation}
  To understand for which $\gamma, j, m$ we have $\varphi_{m, j}(\gamma) \not= 0$ in such an expansion is a difficult  problem (but interesting to my point of view).\\
  The first difficulty is to represent the compact $n$-cycle $\gamma$ (and the horizontal family deduced from it) in order to have a "concrete" integral to compute. Our approach is slightly different and I choose to consider a family of cohomology classes in the fibers $f^{-1}(s)$ induced by a given germ $\omega' \in \Omega^{n+1}_0$, that is to say replace the (holomorphic)  period-integral  $\int_{\gamma(s)} \omega/df$ by the hermitian period 
  $$ \Psi_{\omega, \omega'}(s) := \int_{f^{-1}(s)} \rho \omega/df \wedge \overline{\omega'/df}$$
  where $\rho \in \mathscr{C}^\infty_c(\C^{n+1})$ is identically $1$ near the origin and with a small enough support. Then, thanks to an old result (see \cite{[B.81]}) the function $s \mapsto  \Psi_{\omega, \omega'}(s)$ is smooth on $D^*$ (small enough) and admits an asymptotic expansion in  the space
  $$\vert \Xi\vert^{(n)}_{\mathcal{A}} := \sum_{\alpha \in \mathcal{A}, j \in [0, n]} \C[[s, \bar s]] \vert s\vert^{2\alpha -2}(Log \vert s\vert^2)^j $$
  where $\mathcal{A}$  is a finite subset in $ ]0, 1]\cap \mathbb{Q}$. \\

   To consider  an auxiliary germ $\omega' \in \Omega^{n+1}_0$ is a way to induced an anti-holomorpic family of cohomology classes 
 $\overline{\omega'/df}$ in $H^n(f^{-1}(s), \C)_r$\footnote{The use of the canonical hermitian form on $H^n(F_0, \C)$ is needed for $\alpha = 1$ in the case of an isolated singularity for the eigenvalue $r$ to replace the compact support condition in the usual pairing $H^n(f^{-1}(s), \C)\times H_c^n(f^{-1}(s), \C)  \to \C$ corresponding to the hermitian intersection form for $\alpha \not= 1$  see \cite{[B.85]} or \cite{[B.97]}.}.\\ 
 This presents the advantage to have an integral on the fibers of $f$ (so no more "topological cycle" to find in $f^{-1}(s)$) with also the fact that the non vanishing of some $\varphi_{j,m}(\gamma) $  implies the existence of some $\omega' \in \Omega^{n+1}_0$ such that  the asymptotic expansion of $\Psi_{\omega, \omega'}$ presents a non zero term\footnote{For $\alpha = 1$ we consider only terms with $j \geq 1$, so singular terms.} like $\vert s\vert^{2\alpha - 2}s^m\bar s^{m'}(Log \vert s\vert^2)^j$ for some $m' \in [0, n+1]$.\\

  Using complex Mellin transform of this function (see \cite{[B-M.87]} and \cite{[B-M.89]}) this is equivalent to the existence of a pole of order at least $j$ at the point $-\alpha -m$ for the meromorphic extension of the  function
  \begin{equation}
  F^{\omega, \omega'}_{m'}[\lambda] := \frac{1}{\Gamma(\lambda)} \int_{\C^{n+1}} \vert f\vert^{2\lambda} \bar f^{-m'} \rho \omega \wedge \bar \omega' 
  \end{equation}
  holomorphic for $Re(\lambda) > 0$ and whose meromorphic extension to the complex line is consequence of the existence of the Bernstein-Sato polynomial of $f$ (see \cite{[Be.72]} and  \cite{[B.81]} or \cite{[Bj.93]}).\\

So, our results giving necessary or sufficient conditions for the existence of such a pole in our papers \cite{[B.22]} and \cite{[B.23]} are in fact results  about non zero terms in the asymptotic expansion $(1)$ of the integral period.\\

\section{(a,b)-modules and asymptotic expansions}

Now let me explain the reasons to use (a,b)-modules in this context. The first one comes  from a remark due to Kyoji Saito (in \cite{[KS.83]}, but see also \cite{[Ph.83]} and \cite{[S.89]}) that the Brieskorn module $H^{n+1}_0$ of an isloated singularity germ $f :  (\C^{n+1}, 0) \to (\C, 0)$  has a natural structure of $\C\{\{\partial_s^{-1}\}\}$-module which is free of rank $\mu$ (the Milnor number), where  $\C\{\{\partial_s^{-1}\}\}$ is the $\C$-algebra of Gevrey series\footnote{microdifferential operators of order $\leq 0$ with constant coefficients.}.\\
This structure, with also the obvious $\C\{s\}$-module structure gives a way to look at the Gauss-Manin connexion, that is to say to consider the germ of the differential system at the origin in $\C$ whose multivalued solutions are the period-integrals $(1)$ for any germ $\omega \in \Omega_0^{n+1}$. Recall that in the (totally) isolated singularity case the Brieskorn module\footnote{For any holomorphic germ we may define for each degree $p \in [0, n]$ an (a,b)-module corresponding to the Gauss-Manin system in degree $p$ via the $(p+1)$-th cohomology sheaf of the de Rham complex $(Ker(\wedge df)^\bullet, d^\bullet)$ on $Y = f^{-1}(0)$. See \cite{[B.I]}, \cite{[B.II]} or  \cite{[B-S.04]}. But we discuss here only the case $p = n$.} is simply given by
 $$H^{n+1}_0 := \Omega^{n+1}_0\big/df\wedge d\Omega^{n-1}_0.$$
 In our setting (isolated singularity for the eigenvalue $r := exp(-2i\pi\alpha)$) the Brieskorn module is defined as the quotient by its $b := \partial_s^{-1}$-torsion of the quotient
 $$ \mathcal{H}^{n+1}_0 := \Omega^{n+1}_0\big/ d(Ker(\wedge df)^n \quad {\rm where} \quad \wedge df : \Omega^n_0 \to \Omega^{n+1}_0 .$$
 A serious reason to introduce the (a,b)-module structure (which is nothing else that the simultaneous use of a $\C\{s\}$ and a $\C\{\{\partial_s^{-1}\}\}$-module structures which are "compatible") is the fact that the Gauss-Manin-connection (so the action of $\partial_s$)  is meromorphic and then compels to introduce denominators in $s$ while  the action of $\partial_s^{-1}$ (the primitive  without constant in $s$) is already well defined on the $\C\{s\}$-module $H^{n+1}_0$ and do not need  any denominator in $s$. \\
 But what is the benefice not to invert $s$ and $\partial_s^{-1}$ ??\\
 The answer is very simple : if you look at the action of $s$ and of $\partial_s^{-1}$ on polynomials $\C[s]$ or on holomorphic germs $\C\{s\}$ you immediately see that the filtration by the 
 $b^p\C\{s\}, p \in \mathbb{N}$, defines a natural  filtration  by sub-modules for both actions of $a$ and $b$. It is an easy consequence of  the commutation relation 
 $$ ab - ba = b^2 \quad {\rm where} \ a := \times s \quad {\rm and} \ b = \partial_s^{-1} := \int_0^s $$
 But here it is important to remark that {\bf the filtration defined by $a^pH^{n+1}_0, p \in \mathbb{N}$, is not stable by the action of $b$} in general\footnote{Remark that  the Nullstellensatz gives an integer $N$ such that $f^N\Omega^{n+1}_0 \subset df \wedge \Omega^n_0$ and so  $a^N H^{n+1}_0 \subset b^NH^{n+1}_0$. So the $b$-filtration is finer that the $a$-filtration.}. So the $b$-module structure has to be considered as the "primary one" which is not the usual point of view!\\ 
 Moreover it has been shown by Varchenkho (see \cite{[Va.82]}) that in the case of an isolated singularity this $b$-filtration corresponds to the Hodge filtration  of the mixte Hodge structure defined by Stennbrink (see \cite{[St.76]}) on the cohomology $H^n(F_0, \C)$ in this context.\\
 
 Since the Gauss-Manin connection is regular, every formal solution converges, I decided to explore the role of this bi-structure module in a formal setting and I choose the following definition\footnote{Let me say that the "strange choice" to call $a := \times s$ and $b := \partial_s^{-1} $ was motivated by the following considerations: if you keep the usual notations you risk, at some point, to use an obvious "usual formula" in a context where it is not "obvious". And my interest was to see precisely what comes only from the commutation relation $ab - ba = b^2$ in the study of period-integrals. This choice is also useful to understand where the use of the inverse of $b$ is interesting. It also help me to convince myself that $b$ is the "primary" variable, which is not so easy to accept !}
 \begin{defn}
 An (a,b)-module $E$ is a free finite type $\C[[b]]$-module endowed with a $\C$-linear endomorphism $"a"$ satisfying the commutation relation $ab - ba = b^2$ and which is continuous for the topology defined by $b$-adic filtration $(b^pE)_{p \geq 0}$ of $E$.\\
  In particular we have $aS(b)e = S(b)ae + b^2S'(b)e$ for any $e \in E$ and any $S \in \C[[b]]$.
 \end{defn}
 Then it is very simple to construct an (a,b)-module: take a rank $k$ free $\C[[b]]$-module $E$ with basis $e_1, \dots, e_k$  and define $"a"$ by prescribing the value of $a(e_j) \in E, j \in [1, k]$ as you want. Then this define uniquely an (a,b)-module structure on $E$.\\
 The rank of $E$ on $\C[[b]]$ is called the rank of $E$.\\
 
 Since the most easy case for a regular connection is the case of a  simple pole  (this happens for the Gauss-Manin connection of a quasi-homogeneous polynomial) this leads to the  following definition.
 \begin{defn} 
 An (a,b)-module has a {\bf simple pole} if and only if $aE \subset bE$.
 \end{defn} 
 Remark that under this hypothesis (in fact the existence of an integer $N$ such that $a^NE \subset bE$ is enough) the (a,b)-module $E$ is $a$-complete, so is a module over $\C[[a]]$.\\
 
 Then the following question makes sense and is, of course, fundamental:
 \begin{itemize}
 \item Given a formal simple pole differential system $z\partial_z F(z) = A(z) F(z)$ where $A$ is a $(k, k)$-matrix with entries in $\C[[z]]$ and $F$ an "unknown element"  in $\C[[z]]^k$, does there exists an (a,b)-module $E$ representing such a differential system (in a sense to be precise) ?
 \end{itemize}
 The answer to this question is yes, and the relation from the matrix $A$ to the corresponding (a,b)-module is conceptually simple : take $e := \partial_zF$ and look for $X $ a $(k,k)$-matrix with entries in $\C[[b]]$ such that $ae = bX(b)e$. Proposition 1.1 in \cite{[B.93]}  gives the existence and  uniqueness of  a matrix $X$ such that,  in the formal completion "in $a$" of $E$ (but remember that $E$ is in fact already $a$-complete!)  the equality $ ae = A(a)be$ holds true in $E$. \\
 But to compute explicitely the matrix $X$ from the matrix $A$ is very complicated in general. In some sense, this complexity may explain  partly why the (a,b)-module point of view can really give interesting results and is not a simple game playing with a  change of notations.\\
 
 \parag{Notation} An (a,b)-module is in fact a left module over a non commutative algebra $\mathcal{A}$ which contains both $\C[[b]]$ and $\C[a]$ which is defined as follows:
 $$ \mathcal{A} := \{ \sum_{q = 0}^\infty  P_q(a)b^q\} \quad {\rm where} \quad P_q \in \C[a] \quad \forall q \in \mathbb{N} $$
 where the commutation relation $ab - ba = b^2$ allows to write any element  in $\mathcal{A}$  (in a unique way) as  $\sum_{q=0}^\infty  b^q\Pi_q(a)$ where $\Pi_q$ is also in $\C[a]$. Then it is clear that a left $\mathcal{A}$-module is both a $\C[[b]]$-module and a $\C[a]$-module.\\
 Now an  (a,b)-module is simply  a left $\mathcal{A}$-modules which is free and of finite type on the sub-algebra $\C[[b]]$ of $\mathcal{A}$.\\
  The regularity condition defined below implies that the left action of $\C[a]$ on a regular (a,b)-module extends to  an action of $\C[[a]]$.\\
 
 Then it is easy to define a "regular" (a,b)-module to encode regularity of a connection (recall that Gauss-Manin connection is always regular) as follows:
 \begin{defn}\label{regular}
 An (a,b)-module $E$ is {\bf regular} when it is a sub-module of a simple pole (a,b)-module.
 \end{defn}
 It is easy to see that  when $E$ is regular there exists a smallest canonical simple pole (a,b)-module $E^\sharp$ containing a given regular (a,b)-module $E$  which is the saturation by $b^{-1}a$ of $E$ in $K\otimes_{\C[[b]]} E$ where $K := \C[[b]][b^{-1}]$.\\
 The next step is to define the Bernstein polynomial of a regular (a,b)-module.
 \begin{defn}\label{Bernstein}
 For a simple pole (a,b)-module $E$ the {\bf Bernstein polynomial} of $E$  is the minimal polynomial of the action of $-b^{-1}a$ on the finite dimensional vector space $E/bE$.\\
 For a regular (a,b)-module, its Bernstein polynomial is, by definition, the Bernstein polynomial of $E^{\sharp}$.
 \end{defn}
 
 \parag{Example} Any regular rank $1$ (a,b)-module has always a simple pole and is isomorphic to   $E_\lambda := \C[[b]]e_\lambda$ for some $\lambda \in \C$  where $ae_\lambda = \lambda be_\lambda$. \\
 
 Now to encode the famous result of Kashiwara (see \cite{[K.76]}) on the rationality and  negativity of the roots of the (usual) Bernstein polynomial (see \cite{[Be.72]} or \cite{[Bj.93]})  of a holomorphic germ we define "geometric" (a,b)-modules as follows.
 \begin{defn}
 A regular  (a,b)-module is {\bf geometric} if and only if the roots of its Bernstein polynomial are negative rational numbers.
 \end{defn}
 
 Now I shall show that with such a definition we are quite near from  our initial point of view which was the study of the asymptotic expansion of a period-integral.
 
 \begin{defn}\label{asympt. 0}
 Let $\alpha $ be in $ ]0, 1] \cap  \mathbb{Q}$ and define for $N \in \mathbb{N}$ the (a,b)-module :
 \begin{align*}
 & \Xi_\alpha^{(N)} := \oplus _{j=0}^N \C[[b]]e_j  \quad {\rm where} \\
 & ae_j =  \alpha b(e_j + e_{j-1}) \quad {\rm with \ the \ convention} \quad e_{-1} = 0.
 \end{align*}
 \end{defn}
 It is easy to see that for $\alpha \not= 1$ we may identify $\Xi_\alpha^{(N)}$ with the usual asymptotic expansion $\C[[s]]$-module
 $$ \{ \sum_{j\in [0, N], m \in \mathbb{N}} \varphi_{m, j} s^{\alpha + m -1}(Log\, s)^j \  \} $$
 with the standard actions $a = \times s$ and $b := \int_0^s$. Moreover its has a simple pole and its Bernstein polynomial is equal to $(x+\alpha)^{N+1}$.\ So it is a geometric (a,b)-module.\\
 For any finite subset $\mathscr{A} \subset ]0, 1] \cap  \mathbb{Q}$ and any finite dimensional vector space $V$, the space
 $$ \Xi_{\mathscr{A}}^{(N)} \otimes_{\C} V := \big( \oplus_{\alpha \in \mathscr{A}} \Xi_\alpha^{(N)}\big) \otimes_{\C} V $$
 where $a$ acts by $a \otimes id_V$ and $b$ by $b \otimes id_V$ is again a geometric (a,b)-module. Then we have the following "embedding Theorem"
 
 \begin{thm}\label{asympt. 1}
 Any geometric (a,b)-module $E$ admits an embedding as a sub-module of some $ \Xi_{\mathscr{A}}^{(N)} \otimes_{\C} V$. Moreover we may choose $\mathscr{A}$ as the image in $\mathbb{Q}/\mathbb{Z} \simeq ]0, 1] \cap \mathbb{Q}$ of the opposite of the roots of the Bernstein polynomial of $E$, and $V$ of dimension equal to the rank of $E$\footnote{In fact we have a better estimate for the minimal dimension  of $V$ which is the rank of $S_1(E)$ the semi-simple part of $E$ (see below for the definition).}.\\
 Of course, the converse holds true : any left $\mathcal{A}$ sub-module of some  $ \Xi_{\mathscr{A}}^{(N)} \otimes_{\C} V$ is a geometric (a,b)-module.
 \end{thm}
 
 \section{The use of frescos}
 
 To go now to the study of asymptotic expansions for the periods integrals associated to one germ $\omega \in \Omega^{n+1}_0$ (either for any choice of $\gamma$ or for one choice of $\gamma$ respectively)  the following definition will be used.
 
 \begin{defn}\label{fresco 0}
 Any  sub-module of some $ \Xi_{\mathscr{A}}^{(N)} \otimes_{\C} V$ generated (as a left $\mathcal{A}$-module) by a unique element is called a {\bf fresco}.\\
 When, moreover, the vector space $V$ is one dimensional such a fresco will be called a {\bf theme}.
 \end{defn}
 
 In the setting described at the beginning, the notion of fresco corresponds to the differential  system associated to a given class  $[\omega] \in H^{n+1}_0$. To such a class the associated  fresco 
 $\mathcal{F}_{\omega}$  is the sub-(a,b)-module of the Brieskorn module $H^{n+1}_0$ generated by the class $[\omega]$ (so $\mathcal{F}_\omega := \mathcal{A}[\omega] \subset H^{n+1}_0$).\\
 The notion of theme corresponds to the choice of the class $[\omega]$ with a choice of a $n$-cycle $\gamma$ in $H^n(F_0, \C)$.\\
 
 The frescos have very nice properties :
 \begin{enumerate}
 \item The Bernstein polynomial of a fresco $F$ is the {\it characteristic polynomial} of $-b^{-1}a$ acting on $F^\sharp/bF^{\sharp}$.
 \item Assume that $0 \to E_1 \to E_2 \to E_3 \to 0 $ is an exact sequence of (a,b)-modules\footnote{Note that our definition of an (a,b)-module implies that in such an exact sequence this implies that $E_1$ is {\em normal} in $E_2$ which means that $bE_2 \cap E_1 = bE_1$ in order to avoid $b$-torsion in $E_3$.}. Then if $E_2$ is a fresco,  $E_1$ and $E_3$ are frescos.
 \item In the situation above assuming that $E_2$ is a fresco, the Bernstein polynomials of $E_1, E_2, E_3$ satisfy
 $$ B_{E_2}(x) = B_{E_1}(x + q)B_{E_3}(x) \quad {\rm where} \quad q := rk(E_3).$$
 \end{enumerate}

 The following structure theorem for frescos is very useful\footnote{It allows to give "normal form" for a given fresco and allows  to construct  a versal family (with finitely many parameters) for  themes with given rank and fixed Bernstein polynomial. See \cite{[B.10]}.}
 
 \begin{thm}\label{fresco 1}(see \cite{[B.09]})
 Let $\mathcal{F}$ be a fresco with rank $k$. Then there exists an element in $\mathcal{A}$
  $$P := (a-\lambda_1b)S_1(a - \lambda_2b)S_2 \dots S_{k-1}(a - \lambda_kb) $$
  where $S_1, \dots, S_{k-1}$ are invertible elements in $\C[[b]]$, such that  $\mathcal{F}$ is isomorphic to the quotient $\mathcal{A} \big/ P \mathcal{A}$ as a left $\mathcal{A}$-module. Moreover, in this situation the numbers $-(\lambda_j + j - k)$ are negative rational numbers and the Bernstein polynomial of $\mathcal{F}$ is equal to
  $$B_{\mathcal{F}}(x) = \prod_{j=1}^k (x + \lambda_j + j - k).$$
 \end{thm}

 May be it is time to explain to the reader the kind of result which can be proved using these tools.
 
 \begin{thm}\label{1}(see \cite{[B.22]} and \cite{[B.23]})
 Assume that in the situation described at the beginning the fresco $\mathcal{F}_{\omega}$ associated to $\omega \in \Omega^{n+1}_0$ has a root of order $q \geq 1$  in $-\alpha - \mathbb{N}$. Then there exists a form $\omega' \in \Omega^{n+1}_0$ such that the hermitian period
  $$ \Psi_{\omega, \omega'}(s) = \int_{f^{-1}(s)} \rho\omega/df\wedge \overline{\omega'/df}$$
   has a non zero term in of the form $\vert s\vert^{2\alpha- 2}s^m\bar s^{m'}(Log\vert s\vert^2)^j$ with $j \geq q-1$ if $\alpha \not= 1$ and with $j \geq q$ for $\alpha = 1$.\\
   Conversely, if such a form $\omega'$ exists, then the Bernstein polynomial of $\mathcal{F}_\omega$ has at least $q$ roots (distinct or not)  in $-\alpha - \mathbb{N}$.
   \end{thm}
   
   Note that in this situation, this implies that the (usual) reduced Bernstein polynomial $b_{f,0}$ of $f$ at the origin has at least $q$ roots in $-\alpha - \mathbb{N}$ thanks to \cite{[B.81]}.

 \section{Semi-simple filtration and Higher Bernstein polynomials}
 
 The existence of several roots in $-\alpha - \mathbb{N}$ for the  reduced Bernstein polynomial of  the holomorphic germ $f$ at the origin may imply the existence of some nilpotent order  for the eigenvalue $\exp(-2i\pi\alpha)$  of the monodromy of $f$ at the origin in the case of an isolated singularity.  But {\em this is not always the case}. This phenomena is the same for the Bernstein polynomial of the fresco $\mathcal{F}_\omega$ associated to a germ $\omega \in \Omega^{n+1}_0$.\\
  This explains why in the  theorem above the hypothesis in the direct part (existence of a root with multiplicity $q$...) is stronger than the conclusion of the converse (existence of q roots distinct or not ...).\\
 The goal is now to understand when two different roots in $-\alpha - \mathbb{N}$ for the Bernstein polynomial of the fresco $\mathcal{F}_\omega$ really produces a non zero $Log$-term in the expansion of the 
 period integral of $\omega$. This is, of course,  an attempt to understand the same question for the  roots of the (reduced) usual Bernstein polynomial of $f$ itself. \\
 This problem is clearly linked with the study of the nilpotent part of the action of the {\em logarithm of the monodromy}, which is "naturally" defined  on the saturation $E^\sharp$ of a geometric (a,b)-module via the action of $b^{-1}a$ (see section 6 in \cite{[B.23]}). So this leads to study the semi-simple filtration of a geometric (a,b)-module which we define as follows.
 
 \begin{defn}\label{ss 1}
 An (a,b)-module $E$ is semi-simple if it is a sub-module of a finite direct sum \ $\oplus_{j=1}^N  E_{\lambda_j}$ where $E_\lambda := \C[[b]]e_\lambda$ for $\lambda \in \C$,  the action of $a$ being defined by $ae_\lambda = \lambda be_\lambda$.
 \end{defn}
 
 \begin{prop}\label{ss2}
 For any regular (a,b)-module $E$  there exists a canonical strictly  increasing filtration by normal\footnote{recall that $F \subset E$ is normal when $F \cap bE = bF$. Then the quotient $E/F$ is an (a,b)-module.} sub-modules $(S_j)_{j \in [0, d]}$ with $S_0(E) := \{0\}$ and  also $S_d(E) := E$,  such that the quotients 
  $S_{j+1}(E)\big/S_j(E)$ are semi-simple for each $j $ in  $[0, d-1]$.  The integer $d$ is called the {\bf nilpotent order} of $E$.
 \end{prop}
 
 Since the problem we have in mind concerns only the eigenvalue $exp(-2i\pi \alpha) $ we shall use the following   "decomposition" result associated to the following definition.
 
 \begin{defn}\label{primitive 1}
 A regular (a,b)-module is $\mathscr{A}$-{\bf primitive} if each root of its Bernstein polynomial is in $-\mathscr{A}$ for any subset $\mathscr{A}$ in $\C/\mathbb{Z}$.
 \end{defn}
 
 \begin{prop}\label{primitive 2}
 For any regular (a,b)-module $E$ and any $\mathscr{A} \subset \C/ \mathbb{Z}$ there exists a unique normal sub-module $E_{[\not= \mathscr{A}]}$ in $E$ such that the quotient
 $$ E^{[\mathscr{A}]} := E\big/E_{[\not= \mathscr{A}]} $$
 is $[\mathscr{A}]$-primitive and such that any $[\mathscr{A}]$-primitive sub-module of $E$ does not intersect $E_{[\not= \mathscr{A}]}$. 
 \end{prop}
 
 It is a key point to remark that the semi-simple filtration of a regular (a,b)-module is compatible with its $[\alpha]$-primitive "decomposition" given by the following proposition (see \cite{[B.23]} Lemma 4.2.5)
 
 \begin{prop}\label{HB 0}
 Let $E$ be a regular (a,b)-module and $\mathscr{A}$ be  a subset in $\C/\mathbb{Z}$. Then for each integer  $j$ we have
 $$ S_j(E_{[\mathscr{A}]}) = S_j(E)_{[\mathscr{A}]} \quad {\rm and \ also} \quad S_j(E^{[\mathscr{A}]}) = S_j(E)^{[\mathscr{A}]}.$$
 \end{prop}
 
 Note that in the case where $E$ is a fresco $E^{[\alpha]}$ is also a fresco and  the Bernstein polynomial of $E^{[\alpha]}$ is exactly the $\alpha$-part of the Bernstein polynomial of $E$.\\
 
 Now we can give a "nilpotent weight" to the roots of the Bernstein polynomial of a fresco with the following definition.
 
 \begin{defn}\label{HB 1}
 Let $\mathcal{F}$ be an $\alpha$-primitive fresco with nilpotent order $d$ and for $j \in [1, d]$  and let $\delta_j$ be the rank $\mathcal{F}/S_j(\mathcal{F})$. Define for $j \in [1, d]$  the {\bf $j$-th Bernstein polynomial $B_j(\mathcal{F})$ of $\mathcal{F}$} by the formula
 $$ B_j(\mathcal{F})[x] := \tilde{B}_j(x - \delta_j)$$
 where $\tilde{B}_j $ is the Bernstein polynomial of  $S_j(\mathcal{F})/S_{j-1}(\mathcal{F})$.\\
 Then, for any fresco $\mathcal{F}$ defined its $j$-th Bernstein polynomial by
 $$ B_j(x) := \prod_{\alpha} B_j(\mathcal{F}^{[\alpha]})(x) $$
 where $B_j = 1$ when $\mathcal{F}^{[\alpha]} = \{0\}$ or when $j > d(\mathcal{F}^{[\alpha]})$, the nilpotent order of $\mathcal{F}^{[\alpha]}$.
  \end{defn}
  
  The link with the  Bernstein polynomial in the case of a fresco is given by the following result (\cite{[B.23]} Th. 7.3.2)
  
  \begin{thm}\label{HB 2}
  Let $\mathcal{F}$ be a fresco. The Bernstein polynomial of the fresco $\mathcal{F}$ is given by 
  $$ B_\mathcal{F} = B_1(\mathcal{F})\dots B_d(\mathcal{F})$$
   where $d$ is the nilpotent order of $\mathcal{F}$.\\
Moreover,  each root of each $B_j(\mathcal{F})$ is simple. The degrees of the polynomials $B_j$ are non increasing and are equal to the ranks of the quotients $S_j(\mathcal{F})/S_{j-1}(\mathcal{F})$.
\end{thm}

 Now the main result of our study is the following theorem.
 
 \begin{thm}\label{HB 2}
 In the situation described at the beginning, assume that for some $\omega \in \Omega^{n+1}_0$ the fresco $\mathcal{F}_\omega^{[\alpha]}$ has nilpotent order $d \geq 1$ and that $-\alpha - m$ is  a root of the $d$-th Bernstein polynomial of $\mathcal{F}_\omega^{[\alpha]}$. Then there exists $\omega' \in \Omega^{n+1}_0$ and an integer $m' \in [0, n+1]$ such that the asymptotic expansion of the hermitian period
 $$ \Psi_{\omega, \omega'}(s) := \int_{f^{-1}(s)} \rho \omega/df \wedge\overline{\omega'/df} $$
 has a non zero term of the type \  $\vert s\vert^{2\alpha -2} s^m{\bar s}^{m'}(Log\vert s\vert^2)^{d-1} $  \ for $\alpha \not= 1$ and of the type  \\
  $\vert s\vert^{2\alpha -2}s^m\bar s^{m'}(Log\vert s\vert^2)^d $ for $\alpha = 1$.\\
 Conversely, if such $\omega'$ and $m'$ exists for a given $\omega$, then $\mathcal{F}_\omega^{[\alpha]}$ has nilpotent order at least equal to $d$ and  $d$-th Bernstein polynomial of $\mathcal{F}_\omega^{[\alpha]}$ has a root equal to $-\alpha -m_1$ with $m_1 \in [0, m]$.
 \end{thm}
 
 The reader which is interested with more consequences, more details and proofs may consult \cite{[B.22]} and \cite{[B.23]} where he will also find some examples of computations of the Bernstein polynomial of frescos in explicite examples.\\
 
 To conclude this presentation I must confess that I was not able to prove the last  result above using only the "formal" (a,b)-module point of view presented here. So I was compelled to develop the   {\bf convergent theory} of  (a,b)-modules which does not reduce to an easy exercise on inequalities; for instance see the "division theorem" (in section 2.5 in \cite{[B.23]} ) in the algebra $\mathcal{A}_{conv}$, the convergent analog of $\mathcal{A}$, where the geometric condition is useful and this result  allows to prove the structure theorem for convergent frescos, analog of Theorem \ref{fresco 1} above.

 \end{document}